\documentclass[12pt]{article}
\usepackage{natbib}
\usepackage{fullpage,graphicx,psfrag,theorem}
\usepackage{amsmath,amssymb}
\newtheorem{example}{Example}
\newtheorem{remark}{Remark}
\newtheorem{assumption}{Assumption}
\newtheorem{theorem}{Theorem}

\author{N. Cohen, E. Greenshtein, and Y. Ritov}

\begin{document}
\title { Compound decision in the presence of proxies with an application to spatio-temporal data.}
\maketitle

\def\bb#1{\mathbb{#1}}
\def\bbp{\bb{P}}
\def\bbf{\bb{F}}
\def\en{\infty}
\def\summ#1#2#3{\sum_{#1=#2}^{#3}}
\def\mby{\boldsymbol{y}}
\def\cov{\mathop{\rm cov}}
\newcommand{\bdf}{\boldsymbol}
\begin{abstract}

We study the problem of incorporating covariates in a compound decision setup. It is desired  to estimate the means of $n$  response variables, which are independent and normally distributed, and each is accompanied by a vector of covariates.  We suggest a method that involves non-parametric empirical Bayes techniques and may be viewed as a generalization of the celebrated Fay-Herriot (1979) method.

Some optimality properties of our method are proved. We also compare it numerically with
Fay-Herriot  and other methods,
using a `semi-real' data set that involves spatio-temporal covariates, where the
goal is to estimate certain proportions in many small areas (Statistical-Areas).

\end{abstract}

\section{Introduction}

The main purpose of this paper is to study and demonstrate  how to incorporate compound decision techniques (CD), or almost equivalently, empirical Bayes (EB) methods,
in the presence of explanatory variables.  The ideas of CD/EB were developed in the 1950's by Robbins
(1951, 1955, 1964),   see the review papers by Copas (1969) and Zhang (2003).
Compound decision (or Empirical Bayes)  procedures, were shown to produce very efficient estimators in the simple setup where we have independent observations,
$Y_1,\dots ,Y_n$, $Y_i \sim F_{\mu_i}$,  and it is desired to estimate $\mu_i$, $i=1,\dots ,n$. A major case,
on which we will concentrate, is when $F_{\mu_i}=N(\mu_i,1)$.

We will focus on two types of EB procedures.
One type is Parametric Empirical Bayes (PEB) procedure, where $\mu_i, \; i=1,\dots ,n$ are
assumed to be realizations of independent random variables $M_i, \; i=1,\dots ,n$, $M_i \sim G$, $G = N(0,\tau^2)$, where
$\tau^2$ is unknown and should be estimated from the data.
When $n$ is large, the corresponding estimator
(note, the exact variant of the corresponding estimator, depends on the method of estimating $\tau^2$),
resembles the James-Stein estimator, see e.g., Efron and Morris (1973). The other type is the Non-Parametric
Empirical Bayes (NPEB) procedure, where the above distribution $G$ is a  member of a large non parametric family ${\cal G}$ of
distributions. Two recent NPEB methods and approaches are Brown and Greenshtein (2009) and
Jiang and Zhang (2009).

The advantage of EB procedures, relative to more elementary (e.g., mle) procedures, occures as $n$ grows, and
may become  very significant in
high dimensional problems when $n$ is large (e.g., $n \geq 3$ is needed for ``Stein's paradox'' to hold). A special advantage of NPEB procedures is expected in situations where
the vector $\bdf{\mu}=(\mu_1,\dots ,\mu_n)'$ is sparse, see e.g.,
Greenshtein, Park and Ritov (2008), Brown and Greenshtein (2009).

Since modern statistical problems often involve high dimensional and sparse estimation problems,
EB techniques should be embraced for such purposes, see, e.g., Efron (2003).
However, apart from literature in small area estimation,
e.g., Rao (2003), which follows the seminal paper of Fay and Herriot (1979),
EB is hardly used in modern data analysis. One reason is that in most applied problems, we have
explanatory variables $X_{i1},\dots ,X_{ip}$ for each observation $Y_i$ and in such cases EB has no appeal,
since simple symmetric decision procedures have no appeal. We elaborate in the following.

In our motivating example $Y_i \sim B(m_i,p_i)$, the binomial distribution,  and we need to estimate $p_1,\dots ,p_n$,
certain proportions, in $n$ (small) areas. The values of $p_1,\dots,p_n$ are unknown constants to be estimated. In addition to the sample $Y_1,\dots ,Y_n$, we have a set of variables $\bf{X}_1,\dots ,\bf{X}_n$ (fixed or random, but independent of
$Y_1,\dots ,Y_n$ ) and hopefully $\bf{X}_i$ can serve as proxies to $p_i$, $i=1,\dots ,n$. For example one dimensional
covariates $X_i \sim B(k_i,\tilde{p}_i)$ where $\tilde{p}_i$ are ``typically'' close to $p_i$; alternatively
$\bdf{X}_i$ may be a vector of known parameters of area $i$, that might be ``relevant'' to the parameter of interest
$p_i$, for example the socio-economic level of the region, its size, or mean age. We emphasize two elements.
First, because of the proxies, $Y_1,\dots,Y_n$ cannot be considered as ``permutation invariants'' or ``exchangeable''. Second, we do not believe that the observations follow standard regression models. The covariates are considered as proxies to $p_i$, but they are statistically independent of the $Y$'s (whose only stochastic aspect comes from the binomial sampling), and may be only a rough approximation to $p_1,\dots,p_n$.

{\it Simple symmetric and permutation invariant procedures.} In cases of total ignorance regarding
the parameters of the variables in relation to their identity, e.g., a situation where $ Y_i \sim N(\mu_i,1)$ and there is an exchangeable
multivariate
prior on $(\mu_1,\dots ,\mu_n)$, procedures which are permutation invariant
have a special appeal. Permutation invariant procedures $\Delta$ are such that for every permutation $\pi$,
$$ \Delta(Y_1,\dots ,Y_n)=(a_1,\dots ,a_n) \; \iff  \; \Delta(Y_{\pi(1)},\dots,Y_{\pi(n)})=(a_{\pi(1)},\dots ,a_{\pi(n)});$$
here $a_i \in A$, where $A$ is the ( abstract) action space.
A simple class of exchangeable priors is where $\mu_i$ are realizations of i.i.d $M_i \sim G, \; i=1,\dots ,n$. The optimal procedures then belong to the class of `simple symmetric decision functions', i.e., procedures $\Delta$ which
are of the form: $$\Delta(Y_1,\dots ,Y_n)=(\delta(Y_1),\dots ,\delta(Y_n)),$$ for a given $\delta$.
For natural losses, given $G$, the optimal $\delta$ corresponds to the one dimensional Bayes procedure.
On the relation and asymptotic equivalence between the above two classes, see Greenshtein and Ritov (2009).
Given a loss function, consider an `oracle' who knows the values of
$\mu_1,\dots ,\mu_n$, but is required to use a permutation invariant procedure.
EB and CD procedures may be viewed as an attempt to immitate the (unknown) procedure that
an oracle would use. This
is a very natural goal under `total ignorance' or `exchangeability'.

The appeal in using permutation invariant procedures
and consequently EB procedures, is lost when exchangeability is lost, as in cases where
there are explanatory variables. Assume $n=n_1+n_2$ and it is known that the first $n_1$ observations
(say, hormone measurements), were taken
from men, while the last $n_2$ were taken from women. Applying a permutation invariant procedure is equivalent
to ignoring this potentially important information/explanatory-variable. However not all is lost, one may still apply EB procedure seperately on the first $n_1$ observations and on the last $n_2$ observations. The idea is that after
accounting for the explanatory variable
in this trivial manner, we arrive into (two groups of) exchangeable variables and applying EB procedures
seperately on each group
becomes appealing. In a similar manner, we will account for the information in the explanatory variables
and then, after the information from the explanatory variables is accounted for
and the ``accounted observations'' are closer to being exchangeable, we apply an EB procedure.

EB and CD are closely related notions and approaches. Under an EB formulation the parameters $\mu_i$, $i=1,\dots ,n$
are independent realizations from an unknown distribution $G$ and the aim is to approximate the corresponding Bayes rule; under a CD formulation the aim is to approximate the best decision rule within a class of procedures
(e.g., simple-symmetric, permutation invariant), for the given $\bdf{\mu}=(\mu_1,\dots ,\mu_n)'$. In this paper we
emphasize the CD approach. However, we will often use the more familiar EB notion, motivation and terminology.

Applying a (variant of)  PEB method after accounting for the covariates,
is in the spirit of the paper  of Fay and Herriot, as shown in sub-section 2.2; it is  currently the most common practice. Another approach for inference in the presence of explanatory variables is that
of Lindley and Smith (1972), it is a parametric empirical Bayes approach, though different than
that of Fay and Herriot.

In Section 2, we will suggest how EB could naturally be incorporated in problems with explanatory variables.
We extend the Fay-Herriot approach and present its PEB and NPEB versions.
We show  asymptotic optimality of NPEB.
In section 3, we demonstrate the application of our suggested methods on a ``semi-real'' data, which is based on
the recent Israeli census. The application involves estimation of certain population's proportions in small areas
(Statistical Areas). The explanatory variables available when estimating the proportion $p_i$ in statistical area $i$,
are `Spatial' and `Temporal', based on historical data, and data from neighbors. We elaborate on
comparing PEB procedures, versus the more recent NPEB procedure, suggested
by Brown and Greenshtein (2009).

Our ideas and techniques are meaningful in a genreral setup where $Y_i \sim F_{\mu_i}$, but will be
presented for the case $F_{\mu_i} \equiv N(\mu_i,1)$, $i=1,\dots ,n$. In fact,
as mentioned we will apply our method for estimating
proportions in the setup
where $Y_i \sim B(m_i,p_i)$, but applying an arcsin transformation will bring us to the normal
setup.

 \section{Collections of estimators induced by affine transformations}

The setup we consider is where we observe vectors $\bdf{V}_i=(Y_i,X_{i1},\dots ,X_{ip})$, $i=1,\dots ,n$,
 where $Y_i \sim N(\mu_i,1)$ are independent
response variables, and $X_{ij}$ are explanatory variables independent of $Y_i$, $i=1,\dots ,n$, $j=1,\dots ,p$. Denote by
$X_{n \times p}$ the matrix of the explanatory variables. Denote, $\bdf{Y}'=(Y_1,\dots ,Y_n)$.
The goal is to find a `good' estimator $\hat{\bdf{\mu}}=\hat{\bdf{\mu}}(\bdf{V}_1,\dots ,\bdf{V}_n)$, under the risk
$$E  || \hat{\bdf{\mu}}- \bdf{\mu}||^2_2.$$

In a nutshell the motivation and approach are as follows. Ideally it could be desired to approximate
the  Bayes procedure,  assuming (at leat tactically) that $(\bdf{V}_i,\mu_i)$, $i=1,\dots ,n$, are independent random vectors sampled from  an unknown distribution ${\Gamma}$ that belongs to a large non-parametric family of distributions ${\cal G}$. Then, the goal is to approximate the Bayes decision
$\delta^*= \mathop{\rm \mathop{\rm argmin}}_\delta E_{ {\Gamma}}||\delta(\bdf{V_i})-{\mu_i}||^2$ by $\hat{\delta}^*$,
and let $\hat{\bdf{\mu}}=(\hat{\delta}^*(\bdf{V_1}),\dots ,\hat{\delta}^*(\bdf{V_n}))$.
However, this goal may be too ambitious for $p+1$ dimensional observations $\bdf{V}_i$  when $n$
is moderate due to the ``curse of dimensionality''.
A possible approach, in the spirit of Lindley and Smith (1972), is to assume that $\Gamma$ belongs to a convenient
parametric family, and this way  the``curse of dimensionality'' and other difficulties are  circumvented.
The approach of Fay and Herriot (1979) may also be interpreted this way.
We, on the other hand, aim for the best permutational invariant estimator
with respect to $Z_1,\dots ,Z_n$, where $Z_i$ are one dimensional random variables which are obtained
by a suitable transformation of $(\bdf{V}_1,\dots ,\bdf{V}_n)$. This transformation is estimated from the data.

\subsection{ preliminaries and definitions}

We start from a general point of view, where initially there are no covariates.
We observe independent $Y_i \sim N(\mu_i,1)$, $i=1,\dots ,n$.
Let $\{ T \}$ be a collection of affine
transformations $T(\bdf{Y})=T_{A,B}(\bdf{Y})=A\bdf{Y} - \bdf{B}$, where $A$ is an orthonormal matrix and $\bdf{B}$ is a  vector.
Then $\bdf{Z}=T(\bdf{Y})$ is distributed as a multivariate normal with mean vector denoted $\bdf{\nu}$,
$\bdf{\nu}=A\bdf{\bdf{\mu}}-\bdf{B}$, and covariance matrix the identity. Let $\Delta=\Delta(\bdf{Y})$ be a fixed
estimator of the vector $\bdf{\mu}$, which is not invariant under the group of affine transformations, i.e.,
$\Delta(T(\bdf{Y})) \neq T(\Delta(\bdf{Y}))$. Then, the pair $\Delta$ and $\{ T \}$ defines a
(non-trivial)   class of
decision functions $\{ \Delta_T \}, \; T \in \{ T\}$, $$ \Delta_T(\bdf{Y})= T^{-1}( \Delta(T(\bdf{Y}) ).$$

Let $$ T^{opt}= \mathop{\rm argmin}_{ T \in \{T\} } E_{\bdf{\mu}} || \Delta_T(\bdf{Y})-\bdf{\mu}||_2^2 \equiv
\mathop{\rm argmin}_{ T \in \{T\} } R(T,\bdf{\mu});$$
here $R(T,\bdf{\mu})$ is implicitly defined.

Our goal is to approximate $T^{opt}$, and then estimate $\bdf{\mu}$ by
an approximation of $ \Delta_{T^{opt}}(\bdf{Y})$.

For every  $T \in \{ T \}$, suppose we have a  good estimator $\hat{R}(T,\bdf{\mu})$ for $R(T,\bdf{\mu})$. Let
$\hat{T}=\mathop{\rm argmin}_{ T \in \{T\} } \hat{R}(T,\bdf{\mu})$. The usual approach, which we will
follow,
is to use the estimator
$\hat{ \bdf{\mu} }= \Delta_{\hat{T}}(\bdf{Y})$. When the class $\{ T \}$ is not too large, we expect only a minor
affect of over-fitting, i.e., $ R( \hat{T},\bdf{\mu}) \approx R(T^{opt},\bdf{\mu})$.

\begin{example}[Wavelet transfrom]
The above formulation describes many standard techniques. In fact any harmonic analysis of the data that starts with transforming the data by a standard transformation (e.g., Fourier transform) follows this outline.   A special case  is when $T(\bdf{Y})=A\bdf{Y}$,  where $A$ is the matrix which transforms
$\bdf{Y}$
to a certain  wavelet representation, then, typically, the mean of the transformed vector is estimated and transformed back,
see Donoho and Johnstone (1994). Suppose that, $\{ T \}= \{ A \}$ is a collection of matrices that correspond to
a collection of wavelet bases/``dictionaries''. The problem of finding the most appropriate basis/transformation, is related to that of basis-pursuit,
see e.g., Chen, et.al. (2001).  The permutational invariant and non-linear
decision functions $\Delta$ in those studies is soft/hard-thresolds, Lasso, etc.
As mentioned, procedures of a special interest for us are parametric and non-parametric EB.
\end{example}

\begin{example} [Regression]
Suppose that in addition to $\bdf{Y}$ there is a fixed (deterministic!) matrix $X\in R^{n \times p}$. Consider the class
of transformations $T(\bdf{Y})=\bdf{Y}-\bdf{B}$, $\bdf{B} \in \{ \bdf{B} \}$, where $ \{ \bdf{B} \}$ is the collection of all vectors of the form $\bdf{B}=X\bdf\beta$, $\bdf{\beta} \in R^p$.
Note, in particular, that our transformations are non-random.

\end{example}


\begin{remark} The formulation for
a random set $\{ T \}$, which is independent of $\bdf{Y}$ is just the same.
In the last example when $X_{n \times p}$ is random, we condition on the explanatory variables and arrive to a conditional inference version of the developement
in the sequel.
From a Bayesian perspective, assuming a joint distribution $\Gamma$ as above, conditional independence of
the random set $\{T\}$ and $\bdf{Y}$, conditional on the covariates, follows when we assume that $\bdf{Y}$ and
$X_{n \times p}$ are independent
conditional on $\bdf{\mu}$.
We will remark later on the case where the random set of transformations is `weakly dependent' on $\bdf{Y}$.
\end{remark}

The following  fact is useful. Let $\bdf Z=\bdf Z(T)=\Delta_T(\bdf Y)$. Then $Z_i\sim N(\nu_i,1)$ where $\bdf{\nu}=\bdf{\nu}(T)=\Delta_T(\bdf \mu)$, and
\begin{equation}
R(T,\bdf{\mu}) =E_{\bdf{\mu}} ||\Delta_T(\bdf{Y})-\bdf{\mu}||^2_2 = E_{\bdf{\nu}} ||\Delta(\bdf{Z})-\bdf{\nu}||^2_2
=R(I,\bdf{\nu}). \label{eqn:inv}
\end{equation}
In the last equality $I$ represents the identity transformation.
When there is no real danger of confusion, the dependence on $T$ is suppressed.
We will use equation (\ref{eqn:inv}) later to  establish an estimator $\hat{R}(T,\bdf{\mu})$
for $R(T,\bdf{\mu})$.

The following general three steps method,  for estimating $\bdf{\mu}$, suggests itself.
\begin{description}
\item[Step I:] For every $T$,
estimate $ {R}(T,\bdf{\mu})$ by  $\hat{R}(T,\bdf{\mu})$.

\item[Step II:] Find  $\hat{T}= \mathop{\rm argmin}_T \hat{R}(T,\bdf{\mu})$.

\item[Step III:] Get the estimator: $\hat{\bdf{\mu}}= \hat{T}^{-1}(\Delta(\hat{T}(\bdf{Y}))) \equiv \Delta_{\hat{T}}(\bdf{Y})$.
\end{description}

We summarize. The idea in this subsection is that by an appropriate affine transformation,
that may depend on explanatory variables, we will arrive to a problem which is `easier'
for the procedure $\Delta$ to handle. For example, by choosing an appropriate wavelet basis
we will arrive to a sparse $\bdf{\nu}$, which, roughly, is easier to handle/estimate  the  sparser it is.
More generally, in a rough sense, good permutaion invariant procedures $\Delta$, ``prefer'' to estimate sparse vectors $\bdf{\nu}$,
hence transforming the original problem to a sparse problem is useful.
Indeed accounting for the explanatory variables  in a `good' way, often brings
us to a correponding sparse $\bdf{\nu}$. Moreover, by accounting for explanatory
variables in a good way through a suitable transformation,
we may obtain   (nearly) exchangeable variables;
 whence,
applying a permutation invariant procedure $\Delta$ on the transformed variables becomes natural and appealing.

\subsection{The case where $\Delta$ is  parametric empirical Bayes and the Fay-Herriot procedure.}

In this subsection we study the case where $\Delta$ is a parametric empirical Bayes that corresponds to the prior
$N(0,\tau^2)$, where $\tau^2$ is unknown. When $\tau^2$ is known the corresponding Bayes estimator for $\mu_i$
is $\hat{\mu}_i= \frac{\tau^2}{\tau^2 + 1} Y_i$, and its risk is $\frac{\tau^2}{\tau^2+1}$. When $\tau^2$ is
unknown, we replace $\tau^2$ by its estimate. For our level of asymptotics all consistent estimators
$\hat{\tau}^2$ induce
equivalent estimators  $\hat{\mu}_i= \frac{\hat{\tau}^2}{\hat{\tau}^2 + 1}Y_i$,
and the corresponding estimators are asymptotically equivalent to James-Stein estimator up to $o(n)$, see Efron and Morris (1973).  By working in this level of asymptotic, the considerations in this subsection
are valid for a wide class
of PEB procedures, corresponding to various consistent methods of estimating $\tau^2$, including the J-S procedure.
In particular,   the risk  in estimating a (deterministic) vector $\bdf{\mu}$ by PEB (or James-stein's) method equals:
$$ \frac {n||\bdf{\mu}||^2_2}{||\bdf{\mu}||^2_2 + n}+ o(n). $$

We now examine our three steps estimation scheme, adapted for parametric Empirical Bayes (or, for a James-Stein
estimator $\Delta$). Note that, for every $T$ and the corresponding $\bdf{\nu}$ and $Z_i$ we have:
$ R({I,\nu})= \frac {n ||\bdf{\nu}||^2_2}{||\bdf{\nu}||^2_2 + n}+ o(n). $
Hence a plausible estimator for $R(T,\bdf{\mu})$ is
\begin{equation}
\hat{R}(T,\bdf{\mu})
=\hat{R}(I, \bdf{\nu})
=\max\Bigl\{0,\frac{n(\sum Z_i^2 -n)}{(\sum Z_i^2 -n)+n}\Bigr\}
=\max\Bigl\{0,\frac{n(\sum Z_i^2 -n)}{\sum Z_i^2 }\Bigr\}
\label{eqn:hatr}
\end{equation}

Our three steps scheme adapted for parametric empirical Bayes $\Delta$ is  the following.
\begin{description}

\item[Step I:] For every $T$
estimate $ {R}(T,\bdf{\mu})$ by  (\ref{eqn:hatr}).

\item[Step II:] Find  $\hat{T}= \mathop{\rm argmin}_T \hat{R}(T,\bdf{\mu})$.

\item[Step III:] Get the estimator: $\hat{\bdf{\mu}}= \hat{T}^{-1}(\Delta(\hat{T}(\bdf{Y}))) \equiv \Delta_{\hat{T}}(\bdf{Y})$.

\end{description}

\begin{remark} In the case where   $\{T \}$ corresponds to $\{\bdf{B}= \bdf X\beta: \; \beta \in R^p \}$,
the optimization step II is trivial. We want to minimize the residuals $\sum Z_i^2$. This is achieved
for $\tilde{\bdf{B}}$ which is the projection of $\bdf Y$ on the span of the columns of $\bdf X$, i.e., for
$\hat{T}(Y)=\bdf{Y}-X\hat{\bdf{\beta}}$, where $\hat{\bdf{\beta}}$ is the ordinary least squares estimator.
Upon realizing the last fact, it is easy to see that  our above suggested method  is the method of
Fay and Herriot.
\end{remark}

\subsection{A nonparametric empirical Bayes $\Delta$}

The statements  and development in this sub-section are for nonparametric empirical Bayes procedure
$\Delta$, as in Brown Greenshtein (2009), see appendix.

Let $Z_i \sim N(\nu_i,1)$ be independent. Denote by ${\cal R}(\bdf{\nu})$,
the Bayes risk that corresponds to  the prior which is defined by the empirical
distribution of $\bdf{\nu}$.   Let $f_{\bdf{\nu}}=\frac{1}{n} \sum \phi(z-\nu_i)$, where $\phi$
is the density of a standard normal distribution.
Then
\begin{equation}
{\cal R}(\bdf{\nu})= 1-\int \frac{(f_{\bdf{\nu}}'(z))^2}{f_{\bdf{\nu}}(z)} dz= 1- E_{\bdf{\nu}} \frac{ (f'_{\bdf{\nu}}(Z))^2 }
{(f_{\bdf{\nu}}(Z))^2},
\label{eqn:Fisher}
\end{equation}
see Bickel and Collins (1983).

The following theorem follows from Brown and Greenshtein (2009). It is stated for a triangular array
set up, in order to cover situations of sparse $\bdf{\nu} \equiv \bdf{\nu}^n$. At stage $n$,
$Y_i \sim N(\mu_i^n,1)$ are independent and
for any corresponding sequence $T^n$, $T^n \in \{ T^n \}$, $Z_i \sim N(\nu_i^n,1)$ are independent,
$i=1,\dots ,n$.

\begin{assumption} For every $\alpha>0$ and  every sequence $T^n$ and the corresponding $\bdf{\nu}^n$ we have \newline $\max_i(\nu_i^n)-\min_i(\nu_i^n)=o(n^\alpha)$.
\end{assumption}
\begin{assumption} For some $\alpha_0>0$, $n^{(1-\alpha_0)} {\cal R}(\bdf{\nu}^n) \rightarrow \infty$ for every $T^n$ and corresponding
$\bdf{\nu}^n$.
\end{assumption}

\begin{theorem}\label{th1} Under  Assumptions 1 and 2, for every sequence $T^n$,
\begin{equation}
R(I,\bdf{\nu}^n)=E_{\bdf{\nu}^n} ||\Delta(\bdf{Z})-\bdf{\nu}^n||^2_2 = (1 + o(1))  n{\cal R}(\bdf{\nu}^n)
\end{equation}
\end{theorem}
As explained in the appendix, the procedure $\Delta$ in Brown and Greenshtein requires a bandwidth $h=h_n$,
which approaches slowly to zero. The rate that implies the result in Theorem \ref{th1} is $ h_n \sqrt{\log(n)} \to \infty$.

Given $Y_i \sim N(\mu_i,1)$, and a transformation $T$, $T \in \{T\}$. Let $Z_i$ be the $i'th$
coordinate of $\bdf{Z}=T(\bdf{Y})$.
The last theorem,  and equations (\ref{eqn:inv})
and (\ref{eqn:Fisher}) suggest  the following
estimator $ \hat{R}(T,\bdf{\mu})$ for $ {R}(T,\bdf{\mu})$,
\begin{equation}
\hat{R}(T, \bdf{\mu})=n - \sum {\Large[}\frac { (\hat{f}'_{\bdf{\nu}}(Z_i))  }{\hat{f}_{\bdf{\nu}}(Z_i) }{\Large]}^2,
\label{eqn:hrt0}
\end{equation}
where the density $f_{\bdf{\nu}}$ and its derivative are estimated, for example, by  appropriate kernel estimates.

Only step I of  our general three steps procedure should be adapted, and replaced by:
\begin{description}
\item[Step I:] For every $T$ and corresponding $\bdf{\nu}=\bdf{\nu}(T)$,
estimate $ {R}(T,\bdf{\mu})$ by  (\ref{eqn:hrt0}).
\end{description}

\begin{remark}  Step II could be computationally very complicated when the set $\{T \}$ is large.
In the case where   $\{T \}$ corresponds to
$\{\bdf{B}= \bdf X\beta: \; \beta \in R^p \}$,
a plausible choice, which is computationally convenient is to use the least-squares residuals for $\hat{T}(\bdf Y)$, as in the PEB case. This choice could be far from optimal as will
be demonstrated in the following Examples \ref{ex3} and \ref{ex4} and in the simulations section.

Note, minimizing $R(I,\bdf{\nu})$ with respect to $\bdf{\nu}=\bdf{\nu}(T)$ is equivalent to finding the ``most favorable''
prior, rather than the more conventional task of finding the least favorable prior.

The above method
is reasonable when the class $\{ T \}$ is not too large (in a VC dimension sense) and the overfit affect
is not significant,
otherwise regularization may be required. Those considerations are beyond the scope of our paper.
\end{remark}

Choosing the least squares residuals, as mentioned in the remarks above, may be very inefficient, since it might cause ``smoothing'' of the empirical distribution and
low values of $(f'_{\tilde{\bdf{\nu}}})^2$, which by \eqref{eqn:Fisher} implies  high risk. This could be
caused, e.g., by transforming a sparse structure into a non-sparse one, as in the following Example \ref{ex3}, or
by transforming a structure with well separated groups into  a mixed structure, as in the  Example \ref{ex4}.

\begin{example}\label{ex3} $Y_i \sim N(1,1)$, $i=1,\dots ,2m$, $2m=n$. Suppose we have only one (useless) explanatory variable
$X_i=1$ if $i\leq m$ and 0 otherwise. Projecting $\bdf{Y}$ on $X$, we get $\tilde{\bdf{B}} \approx (1,\dots ,1,0,\dots ,0)'$
and $\bdf{\nu}=\bdf{\mu}- \tilde{\bdf{B}} \approx (0,\dots ,0,-1,\dots ,-1)'$, which is much worst for empirical Bayes estimation than the original $\bdf\mu$: It is easy
to see that $n{\cal R}(\tilde{\bdf{\nu}})=O(n)$, while $n{\cal R}(\bdf{\mu})=0$. From Theorem \ref{th1} we conclude that as
$n \rightarrow \infty$ the advantage of the latter (trivial) transformation compared to the  least squares residuals in terms of
the risk is $o(n)$ compared to $O(n)$.
\end{example}

\begin{example}\label{ex4}
Let $Y_i \sim N(\mu_i,1)$ are independent, where $\mu_i=\mu_1$ for $i=1,\dots ,m$ and $\mu_i=-\mu_1$ for  $i=m+1,\dots,2m=n$.
Suppose $X_i=(\mu_i +W_i) \sim N(\mu_i,1)$,  independent of $Y_i, \; i=1,\dots ,n$.
Let $\tilde{\bdf{\nu}}=\bdf{\mu}-\tilde{ \bdf{ B}}$ where $\tilde{\bdf{B}}$ is the projection of $Y$ on the (random) vector $\bdf X=(X_1,\dots ,X_n)'$. It easy to  check that $\tilde{\nu}_i \rightarrow { \mu_i}/{(\mu_1^2+1)} -  {\mu_1^2}W_i/{(\mu_1^2+1)}$
as $n \rightarrow \infty$. When $\mu_1 \rightarrow \infty$, the empirical distribution of
$\bdf{ \tilde{\nu}} \equiv\bdf{\nu^n}$ converges
to that of a standard normal. The corresponding  Bayes risk ${\cal R}(\bdf{ \tilde{\nu} }^n)$ converges to 0.5. Obviously
the Bayes risk that corresponds to the trivial transformation, for which $\bdf{\nu}^n=\bdf{\mu}^n$,  converges
to zero.
\end{example}

\subsection{Optimality of NPEB $\Delta$.}

Until this point the treatment was for a concrete procedure $\Delta$ and a class $\{ T \}$ of transformations.
The purpose of this section is to advocate the choice of a non-parametric empirical Bayes $\Delta$, which is
denoted $\Delta_{NP}$.

However, as noted, the optimization step (Step II) in the non-parametric approach may be computationally
intensive, so such dominance result might not be enough to persuade that the non-parametric approach might be a good
alternative to the parametric approach and to the Fay Herriot procedure.
In Theorem \ref{th2} below we show that for {\it every} two sequences $\bdf{\mu}^n$ and  $T^n$,
the sequence of estimators, that is obtained by coupling $T^n$ with $\Delta_{NP}$,
asymptotically dominates the sequence which is obtained when coupling the same $T^n$ with any other
sequence of permutation invariant procedures  $\Delta^n$.

 Given a procedure $\Delta$, a transformation
$T$, and a mean vector $\bdf{\mu}$, the corresponding risk is denoted for simplicity as
$R_{\Delta}(T,\bdf{\mu}) \equiv R(T,\bdf{\mu})$ as before; for the case of nonparametric EB procedure $\Delta_{NP}$,
the corresponding risk is denoted  $R_{NP}(T,\bdf{\mu})$. Similarly to the previous sub-section our asymptotic analysis is  of a triangular array setup.

\begin{theorem}\label{th2}   Let $\bdf{\mu}^n$, $\Delta^n$ and $T^n$ be arbitrary sequences.
Assume that for each $n$ the procedure $\Delta^n$ is simple symmetric. Further assume Assumptions 1,2. Then:
$$ \limsup \frac{ R_{NP}(T^n,\bdf{\mu}^n) }{  R_{ \Delta^n }(T^n,\bdf{\mu}^n)    } \leq 1. $$
\end{theorem}

\noindent\emph{Proof:} Follows from Brown and Greenshtein (2009) and Theorem 1.
Note that, the risk of the optimal simple symmetric procedure equals $n{\cal R}(\bdf{\nu}^n)$.\newline


\noindent{\bf Conjecture:}  It seems that in Theorem 2,
the condition that $\Delta^n$ are simple symmetric for every $n$, may be replaced by the weaker condition, that $\Delta^n$ are permutation invariant for every $n$. This should follow by an
equivalence result in the spirit of Greenshtein and Ritov (2009), though stronger.
Note, the equivalence result in Greenshtein and Ritov (2009) would suffice under the assumption  that $\max_i(\nu_i^n)-\min_i(\nu_i^n)=O(1)$;
however, Assumption 1 allows a higher order.

\subsection{Remark}
The following remark is for the case in which we are mainly interested, where
$\{T \}$ corresponds to $\{ \bdf{B} = X\beta \}$.
 Denote $\bdf{B}=(B_1,\dots ,B_n)'$. In the application we have in mind
 the set $\{ T \}$ may be random since $X_{ij}$ could be random. When the random set of transformations is independent
 of $\bdf{Y}$, our above treatment applies by conditioning on the explanatory variables.  We will be interested
 in situations where the random set $\{ T \}$ may depend on  $\bdf{Y}$, however we will require that
 $Y_i$ is independent of $X_{i1},\dots , X_{ip}$ for each $i$. Then the conditional distribution of
 $Z_i$ conditional on $X_{i1},\dots ,X_{ip}$ is $N(\nu_i,1)$,
 where $(\nu_1,\dots ,\nu_n)'=\bdf{\nu}=A \bdf{\mu} - \bdf{B} $ as before. When the dependence
 of $Y_i$ on $X_{j1},\dots ,X_{jp}$, $j \neq i$ is not too heavy, a natural goal is still to  try to approximate
 the best decision function for estimating $\nu_i$ among the decision functions which are simple symmetric with respect to $Z_1,\dots ,Z_n$. The conditional marginal distribution of $Z_i $, $i=1,\dots ,n$  is still
 $N(\nu_i,1)$ (i.e., the conditional distribution of $Z_i=Y_i -B_i$ conditional upon $(X_{i1},\dots ,X_{ip})$); however, we may not treat them as independent observations. Thus, the rates of estimating $f_{\bdf{\nu}}$
 and its derivative may become slower, and for heavy dependence, Theorems 1 and 2 might not hold. Similarly,
 rates of estimation of $\tau^2_n$,   in order to apply the PEB procedure, could be slow. However, when the dependence is not ``too heavy'' we may expect Theorems 1 and 2 to hold under the assumption that $Y_i$ is independent of $X_{i_1},\dots , X_{ip}$ for each $i$.

\section{Simulation}

\subsection{Preliminaries}

The city of Tel Aviv is divided into 161 small areas called ``statistical areas'',
each area belongs to a sub-quarter that includes about four additional statistical areas. In the recent Israeli census
the following the proportion $p_i$  of people who are registered in area $i$ among those who live in area $i$, , $ i=1,\dots 161$, were of interest as part of  the process of estimating the population in each
statistical area. The estimated $p_i$, $i=1,\dots ,n$ are used in order to adjust the administrative-registration
counts and get population estimates for each area.
We will not elaborate on it. In our simulation we use as  $p_1,\dots,p_n$ their value as estimated in the recent census (where about $20\%$ of the population was sampled). The mean of $p_i, \; i=1,\dots ,161,$ is 0.75 and their standard deviation is 0.13, their histogram is roughly bell shaped.

We will present a simulation study in which $p_i$, $i=1,\dots ,161$ are estimated based on
samples of size $m_i$ and corresponding simulated independent $Y'_i$,
$Y'_i \sim B(m_i,p_i)$. Here $Y'_i$ is the number of people in the sample from area $i$, which are registered
to area $i$.

In addition we will simulate covariates. We will simulate
temporal variables that correspond to historical data from each area $i$, and
spatial covariates, that correspond to samples from the neighboring areas of each area $i$.
In the following we will explore simulations and scenarios for the cases of: only temporal covariates,
only spatial covariates, and both temporal and spatial covariates.
 We will compare the performance of PEB, NPEB and other methods.
In all the analyzed situations, we will simulate binomial observations with sample size $m_i \equiv m$, for  $m=25, 50, 100$.

In order to reduce this setup to the above normal case, we apply an arcsin transformation
on our binomial observations $\tilde Y_i$, $i=1,\dots ,n$,   as
in Brown (2008).
Specifically, let
\begin{equation}
Y_i= \sqrt{4m}\arcsin \Bigl(\sqrt{ \frac{{\tilde Y_i}+ 0.25}{m+0.5} } \; \Bigr).\label{eqn:bin1}
\end{equation}
Then, $Y_i$ are distributed approximately as $N(\sqrt{4m} \arcsin ( \sqrt{p_i} ),1)$.
We estimate $\mu_i=E(Y_i)$,  by $\hat{\mu}_i$, $i=1,\dots ,n$, as explained in sub-sections 2.3 and 2.3, and then let the estimate of $p_i$, $i=1,\dots ,161$ equal,
\begin{equation} \hat{p}_i=   ( \sin(  \frac{ \hat{ \mu} _i}{ \sqrt{4m} })   )^2. \label{eqn:bin2} \end{equation}
Let $\bdf{p}=(p_1,\dots ,p_n)$ and $\bdf{\hat{p}}=(\hat{p}_1,\dots ,\hat{p}_n))$.
We evaluate the performance of an estimator according to the risk
$$ E_{\bdf{p}} ||\bdf{ \hat{p}} - \bdf{p}||^2_2.$$   The risk is approximated through  1000 simulations
for each entry in the tables in the sequel.
A different parametric EB approach for estimating proportions in small areas, that involves a logistic regression model, may be found in Farell, et.al.

\subsection{Temporal Covariates}

We introduce now simulated scenarios with only Temporal  covariates.
We think of a process where each year a sample of size $m$ is taken from each area.
Suppose we use the records of the previous three years as covariates. Let $\tilde T_i$ be the number
of people among the $3m$ that were sampled in the previous  three years from area $i$, which were registered to the
area. Although $\tilde T_i$ might be better modeled as a binomial mixture, we will  model $\tilde T_i$ as
$B(3m,p_{it})$ for simplicity.  In order to (hopefully) have a linear relation between the response and explanatory variable,
we define our temporal covariates as:
\begin{equation}
T_i= \sqrt{4m}\arcsin \Bigl(\sqrt{ \frac{{\tilde T_i}+ 0.25}{3m+0.5} } \; \Bigr).\label{eqn:temp}
\end{equation}
Note, if there is hardly any change from the previous three years to the current year in area $i$, we will
have $p_i \approx p_{it}$ and $E(T_i) \approx E(Y_i)$.

In the following we will simulate two scenarios. One scenario is of no-change where  $p_{it}=p_i$, $i=1,\dots ,161$. The other scenario is of a few abrupt changes; specifically, $p_i=p_{it}, \; i=17,\dots ,161$, however $p_{it}=0.3 < p_i$ for $i=1,\dots ,16$. Such abrupt changes could occur in areas that went  in previous years through a lot of
building, internal immigration and other changes.

Since the empirical distribution of $E(Y_i)$ is roughly bell-shaped it is expected that the PEB method will work well
in the no-change scenario. Under the  few abrupt changes, an advantage of the NPEB procedure will be observed.

As mentioned in Section 2, the optimization step of the NPEB procedure is difficult. We will try two candidate
transformations $Y-\bdf{B}^i$, $i=1,2$, coupled with the NPEB, the corresponding methods are denoted NPEB1 and NPEB2.
NPEB1 corresponds to the least-squares/Fay-Herriot transformation, while NPEB2 corresponds to the transformation
$Z_i=Y_i- T_i$ (i.e., $\bdf{B}^2=(T_1,\dots ,T_n)'$ ).
The later transformation, although still sub-optimal when coupled with a NPEB $\Delta$,
could occasionally perform better than the former, as also indicated  by Examples 3 and 4.
In addition to comparing the risks of the PEB and NPEBi, $i=1,2$ methods, we will also compare the
the risk of the naive estimator, and of the regression estimator. The regression estimator
estimates $\hat{\mu}_i$ through the least squares predictor (i.e., $\hat{\bdf{\mu}}=X \hat{\beta}$), however it does not apply an additional
PEB or NPEB stage. The Naive estimator simply estimates $p_i$ by the corresponding sample proportion.

The no-change scenario is presented in Table 1. Each entry  is based on 1000 simulated realizations.
Under no-change the temporal covariate is very helpful, and even the regression-estimate,
i.e. least squares linear predictor is doing very well. Over all, the Naive estimator is the worst,
NPEB1, NPEB2 and Regression are about the same, while the PEB is moderately better than the other methods.

\begin{table}
\caption{  }\label{Table1}
\vspace{1ex}\begin{center}
\begin{tabular}{|l|rrrrr|}
\hline&&&&&\\
& Naive
& Reg
&NPEB1
& NPEB2
& PEB
\\\hline&&& &&\\
$m=25$ &  1.12
  & 0.33
  & 0.35
  & 0.37
  &  0.27\\
 $m=50$ & 0.56& 0.17 & 0.18 & 0.18 & 0.14\\
$m=100$ & 0.28 &0.092 &  0.093 & 0.093 &0.073
 \\&&& &&\\\hline

\end{tabular}
\end{center}
\end{table}

Next we consider the scenario of a few abrupt changes.
In this scenario the  regression by itself is performing the worst, however an additional
EB step is helpful. Here the NPEB2 procedure is the best, see Table 2.

\begin{table}
\caption{  }\label{Table2}
\vspace{1ex}\begin{center}
\begin{tabular}{|l|rrrrr|}
\hline&&& &&\\
& Naive
& Reg
&NPEB1
& NPEB2
& PEB
\\\hline&&& &&\\
$m=25$ &  1.12
  & 1.66
  & 0.75
  & 0.49
  &  0.68\\
 $m=50$ & 0.56& 1.64 & 0.46 & 0.22 & 0.42\\
$m=100$ & 0.28 &1.62 &  0.26 & 0.11 &0.24
 \\&&& &&\\\hline

\end{tabular}
\end{center}
\end{table}

\subsection{Spatial Covariates}

In this section we simulate a scenario with spatial covariates.
Tel-Aviv is divided into sub-quarters, where a few statistical areas
define a sub-quarter. Each sub-quarter is defined by about 5 statistical areas.
For every $i=1,\dots ,161,$ we define the neighborhood of area $i$'', as all the statistical areas
{\it other than} area $i$, that belong to the same sub-quarter as area $i$.

Based on the census we have  good estimates for
$p_{is}$- the proportion of people living in the neighborhood of area $i$, who are registered
to their areas. Those estimates are treated as the ``real'' values in our simulations.
The correlation between $p_i$ and $p_{is}$, $i=1,\dots ,161$ is 0.62.

For simplicity we will assume that for each $i$, the size of the sample from
the neighborhood of area $i$ is $4m$.
Let $\tilde S_i$ be the number of  people sampled from the neighborhood of $i$, who are registered to their area.
Although $\tilde S_i$ might be better modeled as a binomial mixture, we will  model $\tilde S_i$ as
$\tilde S_i \sim B(4m,p_{is})$ for simplicity. As in the case of Temporal covariates we define the Spatial covariate for area $i$ as: \begin{equation}
S_i= \sqrt{4m}\arcsin \Bigl(\sqrt{ \frac{{\tilde T_i}+ 0.25}{4m+0.5} } \; \Bigr).\label{eqn:spat}
\end{equation}
As in the temporal case we will consider two NPEB estimates, corresponding to the projection/Fay-Herriot and to the $Z_i=Y_i-S_i$ transformations. The results of our simulations are summarized in  Table 3.
The advantage of the EB procedures is more noticeable for small $m=25$. The explanation is the following. Since the temporal
covariate is not very strong, $\bdf{\nu}$-the mean of the transformed variables  is not too sparse.
When $m$ is large, under the
scale which is induced by the variance
of $Z_i$, the points $\nu_i, \; i=1,\dots ,n,$ may be viewed as isolated (i.e., extremely non sparse)  and the smoothing of the EB is hardly effective. Hence the EB methods behave roughly like the Naive estimator.

\begin{table}
\caption{  }\label{Table3}
\vspace{1ex}\begin{center}
\begin{tabular}{|l|rrrrr|}
\hline&&& &&\\
& Naive
& Reg
&NPEB1
& NPEB2
& PEB
\\ \hline&&& &&\\
$m=25$ &  1.12
  & 1.41
  & 0.72
  & 0.75
  &  0.64\\
 $m=50$ & 0.56& 1.34 & 0.44 & 0.44 & 0.40\\
$m=100$ & 0.28 &1.31 &  0.26 & 0.28 &0.23
 \\&&& &&\\\hline
\end{tabular}
\end{center}
\end{table}

One could wonder whether the spatial covariates are helpful at all, for the non parametric empirical Bayes, i.e,
may be it is better not to transform the data at all and to apply $\Delta_{NP}$ on the
original data taking $T=I$ and $\bdf{\nu}=\bdf{\mu}$. However this option is
slightly worst than the above ones. The simulated risks that correspond to $m=25,50, 100$ are 0.84 , 0.5 and 0.28.

\subsection{Spatial and Temporal Covariates.}
In this sub-section we study the performances of our estimators when both the temporal and spatial variables are introduced. As before we will apply the projection transformation for the NPEB estimator. However, we will
also try the transformations $Z_i= Y_i - ( \alpha S_i +(1-\alpha) T_i)$, for $\alpha = 0, 0.3, 0.6$.
The corresponding estimators
are denoted: NPEB1 (for the projection), NPEB2, NPEB3 and NPEB4 correspondingly.
For the temporal covariates we simulate the scenario of 16 abrupt changes, the spatial covariates are as before.
As may be expected, since the spatial covariate is weak  relative to the temporal, accounting for it causes
extra unnecessary smoothing. For the non-parametric EB procedure, indeed  NPEB2 that corresponds to
$\alpha=0$ has the best performance, which is also the optimal among all the seven methods.

\begin{table}
\caption{  }\label{Table4}
\vspace{1ex}\begin{center}
\begin{tabular}{|l|rrr rrrr|}
\hline&&& && &&\\
& Naive
& Reg
&NPEB1
& NPEB2
&NPEB3
&NPEB4
& PEB
\\&&& &&&&\\\hline&&& &&&&\\
$m=25$ &  1.12
  & 1.13
  & 0.65
  & 0.49
  & 0.54
  & 0.55
  & 0.58 \\
 $m=50$ & 0.56& 1.06 & 0.4 & 0.22 & 0.28 & 0.38 & 0.37\\
$m=100$ & 0.28 &1.03 &  0.24 & 0.11 & 0.15 & 0.22 & 0.22 \\
&&& &&&&\\\hline
\end{tabular}
\end{center}
\end{table}

\section{Appendix}

{\it NPEB procedure. } We will follow the approach  of Brown and Greenshtein (2009), see that paper
for further details.

Assume $Z_i \sim N(\nu_i,\sigma^2)$, $i=1,\dots ,n$, where $\nu_i \sim G$.

Let $$ f(z)=  \int \frac{1}{\sigma}\varphi {\Huge ( } \frac{ z- \nu}{\sigma} {\Huge)}dG(\nu).$$
It may be shown that the normal Bayes procedure denoted $\delta^G_N$, satisfies:
\begin{equation} \delta_N^G(z)= z + \sigma^2 \; \frac{f'(z)}{f(z)}.  \label{eqn:NEB} \end{equation}
The procedure studied in Brown and Greenshtein (2009), involves an estimation of $\delta_N^G$,
by replacing $f$ and $f'$ in (\ref{eqn:NEB}) by their kernel estimators  which are derived through
a normal kernel with bandwidth $h$. Denote the kernel estimates by $\hat{f}_h$ and $\hat{f}'_h$
we obtain the decision function, $(Z_1,\dots,Z_n)\times z \mapsto R$:
\begin{equation} {\delta}_{N,h}(z)= z + \sigma^2 \; \frac{\hat{f}_h'(z)}{\hat{f}_h(z)}.  \label{eqn:ENEB} \end{equation}

A suitable (straightforward) truncation is applied when estimating the corresponding mean of points
$Z_i$ for which $\hat{f}(Z_i)$ is too close to zero and consequently $|\delta_{N,h}(Z_i)-Z_i|> 2\log(n)$.
We did not apply such truncation in our simulations in this paper.
The default choice for the bandwidth $h \equiv h_n$, suggested by Brown and Greenshtein is $1/\sqrt{\log(n)}$.
See also, a cross-validation method for choosing $h$, suggested by Brown, et.al., (2010), together with
some suggested improvements of the above procedure. In our numerical study, $n=161$ and we chose $h=0.4$.
The procedure is not too sensitive to the choice of $h$.

\mbox{}\par

\noindent {\bf \Large References}

\begin{list}{}{\setlength{\itemindent}{-1em}}
\item
Bickel, P. J. and Collins, J.R. (1983). Minimizing Fisher information over mixtures of distributions.
{\it Sankhya} {\bf Vol 45}, No. 1, p 1-19.
\item
Brown, L. D. (2008). In-Season Prediction of Bating Averages: A
field test of Simple Empirical Bayes and Bayes Methodologies. {\it
Ann. of App. Stat.} {\bf 2} 113-152.
\item
Brown, L.D. and Greenshtein, E. (2009). Non parametric empirical Bayes and compound decision
approaches to estimation of high dimensional vector of normal means. {\it Ann. Stat.} {\bf 37}, No 4, 1685-1704.
\item
Brown, L.D, Greenshtein, E. and Ritov, Y. (2010). The Poisson compound decision problem revisited. Manuscript.
\item
Chen, S.S, Donoho, D.L., Saunders, M. A (2001). Atomic decomposition by basis pursuit. {\it SIAM Rev.}
{\bf Vol 43}, {\bf Issue 1}, 129-159.
\item
Copas, J.B. (1969). Compound decisions  and empirical Bayes (with discussion). {\it JRSSB} {\bf 31} 397-425.
\item
Donoho, D.L. and Johnstone, I.M (1994). Ideal spatial adaptation by wavelet shrinkage. {\it Biometrika} {\bf 81}
No. 3, 425-455.
\item
Efron, B. and Morris, C. (1973). Stein's estimation rule and its competitors- an Empirical Bayes approach.
{\it JASA} {\bf 68} 117-130.
\item
Efron, B. (2003). Robbins, Empirical Bayes, and Microarrays (invited paper). {\it Ann.Stat} {\bf 31}, No. 2, 364-378.
\item
Fay, R.E. and Herriot, R. (1979). Estimates of income for small places: An application of James-Stein
procedure to census data. {\it JASA}, {\bf 74}, No. 366, 269-277.
\item
Farrell, P.J., MacGibbon, B., Tomberlin, T.J. (1997). Empirical Bayes estimators of small area proportions
in multistage designs. {\it Stat. Sinica} {\bf 7} 1065-1083.
\item
Greenshtein, E., Park, J., Ritov, Y. (2008). Estimating the mean of high valued observations in high dimensions.
{\it Journal of Stat. theory and pract.} {\bf 2} No.3, 407-418.
\item
Greenshtein, E. and Ritov, Y. (2009). Asymptotic efficiency of simple
decisions for the compound decision problem. The 3'rd Lehmann Symposium. IMS Lecture Notes Monograph Series, J.Rojo, editor.
\item
Lindley, D.V. and Smith, A.F.M. (1972). Bayes estimates for the linear model. {\it JRSSB} {\bf 34}, No.1, 1-41.
\item
Rao, J.N.K. (2003). Small area estimation. Wiley \& Sons, New Jersey.

\item
Robbins, H. (1951). Asymptotically subminimax solutions  of compound decision problems. {\it Proc. Second Berkeley Symp.} 131-148.
\item
Robbins, H. (1955). An Empirical Bayes approach to statistics. {\it  Proc. Third Berkeley Symp.} 157-164.
\item
Robbins, H. (1964). The empirical Bayes approach to statistical decision problems. {\it Ann.Math.Stat.}
{\bf 35}, 1-20.
\item
Jiang, W. and Zhang, C.-H. (2009). General maximum likelihood
empirical Bayes estimation of normal means. {\it Ann. Stat.} {\bf 37}, No 4, 1647-1684.
\item
Zhang, C.-H.(2003). Compound decision theory and empirical Bayes methods.(invited paper). {\it Ann. Stat.} {\bf 31} 379-390.

\end{list}

\end{document}